\documentclass{article}
\usepackage{mathrsfs}
\usepackage{amsmath}
\usepackage{amssymb, amsmath}
\usepackage{graphicx}
\catcode`\@=11 \@addtoreset{equation}{section}

\catcode`\@=12
\usepackage{color}

\newtheorem{Theorem}{Theorem}[section]
\newtheorem{Proposition}{Proposition}[section]
\newtheorem{Lemma}{Lemma}[section]
\newtheorem{Corollary}{Corollary}[section]
\newtheorem{Definition}{Definition}[section]

\newtheorem{Remark}{Remark}[section]

\newcommand{\bTheorem}[1]{
\begin{Theorem} \label{T#1} }
\newcommand{\eT}{\end{Theorem}}

\newcommand{\bProposition}[1]{
\begin{Proposition} \label{P#1}}
\newcommand{\eP}{\end{Proposition}}

\newcommand{\bLemma}[1]{
\begin{Lemma} \label{L#1} }
\newcommand{\eL}{\end{Lemma}}

\newcommand{\bCorollary}[1]{
\begin{Corollary} \label{C#1} }
\newcommand{\eC}{\end{Corollary}}

\newcommand{\beq}{\begin{equation}}
\newcommand{\eeq}{\end{equation}}
\newcommand{\bFormula}[1]{
\begin{equation} \label{#1}}
\newcommand{\eF}{\end{equation}}

\newcommand{\f}{\frac}
\newcommand{\Om}{\Omega}

\newcommand{\vr}{\varrho}

\newcommand{\vu}{\vc{u}}

\newcommand{\vc}[1]{{\boldsymbol #1}}
\newcommand{\Div}{{\rm div}}
\newcommand{\Grad}{\nabla}

\newcommand{\dx}{{\rm d} x}

\newcommand{\dt}{{\rm d} t }
\newcommand{\ds}{{\rm d} s}

\newcommand{\dxdt}{\dx\dt}

\font\F=msbm10 scaled 1000
\newcommand{\R}{\mbox{\F R}}

\newcommand\Cbox[2]{%
    \newbox\contentbox%
    \newbox\bkgdbox%
    \setbox\contentbox\hbox to \hsize{%
        \vtop{
            \kern\columnsep
            \hbox to \hsize{%
                \kern\columnsep%
                \advance\hsize by -2\columnsep%
                \setlength{\textwidth}{\hsize}%
                \vbox{
                    \parskip=\baselineskip
                    \parindent=0bp
                    #2
                }%
                \kern\columnsep%
            }%
            \kern\columnsep%
        }%
    }%
    \setbox\bkgdbox\vbox{
        \color{#1}
        \hrule width  \wd\contentbox %
               height \ht\contentbox %
               depth  \dp\contentbox
        \color{black}
    }%
    \wd\bkgdbox=0bp%
    \vbox{\hbox to \hsize{\box\bkgdbox\box\contentbox}}%
    \vskip\baselineskip%
}

\begin{document}


\title{\bf Global weak solutions to the one-dimensional compressible heat-conductive MHD equations without resistivity}

\author{Yang Li \\ Department of Mathematics, \\ Nanjing University, Nanjing 210093, China \\ lynjum@163.com \\
Yongzhong Sun \\ Department of Mathematics, \\ Nanjing University, Nanjing 210093, China \\ sunyz@nju.edu.cn}

\maketitle
{\bf Abstract: }{We investigate the initial-boundary value problem for one-dimensional compressible, heat-conductive, non-resistive MHD equations of viscous, ideal polytropic fluids in the Lagrangian coordinates. The existence and Lipschitz continuous dependence on the initial data of global weak solutions are established. Uniqueness of weak solutions follows as a direct consequence of stability. }  \\


\section{Introduction}
Magnetohydrodynamics (MHD) studies the mutual interaction between moving, conducting fluids and magnetic fields. It has a wide range of applications in liquid metals, plasmas, and geophysics. In accordance with Newton's laws of motion and Maxwell's laws of electrodynamics, a simplified and classical model for three-dimensional compressible MHD equations takes the following form in Eulerian coordinates (see \cite{HC}):
\beq\label{g1}
\vr_t + \Div(\vr\vu)=0,
\eeq
\beq\label{g2}
(\vr\vu)_t+\Div(\vr\vu \otimes \vu)+\Grad{p}=\Div\mathbb{S}
+(\Grad \times \vc{b})\times \vc{b},
\eeq
\beq\label{g3}
(\vr e)_t+\Div(\vr e \vu)+\Div\vc{q}=\mathbb{S}:\nabla\vu-p\Div\vu+\lambda|\nabla\times\vc{b}|^2,
\eeq
\beq\label{g4}
\vc{b}_t=\Grad \times (\vu \times \vc{b})-\lambda \Grad \times (\Grad \times \vc{b}),\,\Div \vc{b}=0.
\eeq
Here the unknown functions $\vr,\vu\in \R^3,p,e,\vc{b}\in \R^3$ denote the density of the fluid, the velocity field, the pressure, the internal energy and the magnetic field, respectively. The nonnegative constant $\lambda$ represents the magnetic diffusion of the field $\vc{b}$.

The symbol $\mathbb{S}$ stands for the viscous stress tensor:
\[
\mathbb{S}=\eta(\Div\vu)\mathbb{I}+\nu\left(\nabla\vu+(\nabla\vu)^\top\right),
\]
where $\nu$ and $\eta$ are the shear viscosity coefficient and the bulk viscosity coefficient respectively, satisfying physical conditions:
\[
\nu>0,3\eta+2\nu\geq0.
\]
The transpose of matrix $\nabla\vu$ is denoted by $(\nabla\vu)^\top$ and $\mathbb{I}$ is the identity matrix of order three.

The heat flux $\vc{q}$ is subject to Fourier's law:
\[
\vc{q}=-\kappa\nabla\vartheta,
\]
where the positive constant $\kappa$ denotes the heat conductivity coefficient and $\vartheta$ is the absolute temperature of the fluid.

There are many mathematical investigations to the Navier-Stokes-Fourier MHD equations (\ref{g1})-(\ref{g4}) with positive resistivity coefficient $\lambda$. Remarkably, based on the pioneering work of P. L. Lions \cite{LS} and E. Feireisl et al. \cite{FNP} on the three-dimensional compressible isentropic Navier-Stokes equations, Hu and Wang \cite{HW} proved the existence of global weak solutions to the three-dimensional compressible MHD equations with finite energy initial data. B. Ducomet and E. Feireisl \cite{DF} obtained the global weak solution to a more physical model with the influence of high-temperature radiation and the gravitational force.

However, the corresponding results are not so fruitful provided the resistivity coefficient $\lambda$ vanishes. Very recently, Wu and Wu \cite{WW} proved the global small solutions to the two-dimensional compressible MHD equations when the initial data is close to an equilibrium state. This is a compressible counterpart compared with incompressible fluids. Specifically, Lin et al. \cite{LXZ} showed the global small solution to the two-dimensional incompressible MHD equations without magnetic diffusion; see also \cite{RWZ,XZ,ZT} for more results. As a matter of fact, vanishing resistivity gives rise to extra difficulty in deriving necessary a priori estimates and obtaining global solutions, due to the loss of regularizing effect and the strong interaction between the fluid motion and the magnetic field. Whether there  exists a global weak solution or not to the system (\ref{g1})-(\ref{g4}) (with $\lambda = 0$) is still unknown.

In this paper, we restrict ourselves to the MHD equations (\ref{g1})-(\ref{g4}) for an ideal, polytropic gas without resistivity in one spatial dimension. Specifically, the constitutive relations for the pressure $p$ and the internal energy $e$ are given by
\[
p=R\vr\vartheta,\,e=c_V\vartheta,
\]
where $R>0$ is the perfect gas constant and $c_V>0$ means the specific heat at constant volume. Furthermore, based on a specific choice of dependent variables:
\[
\vr=\vr(x,t),\,\vu=(u(x,t),0,0),\,\vartheta=\vartheta(x,t),\,\vc{b}=(0,0,b(x,t)),
\]
where $x\in \R$ is the spatial variable and $t\geq0$ denotes the temporal variable, (\ref{g1})-(\ref{g4}) (with $\lambda = 0$) are transformed to
\beq\label{g5}
\vr_t+(\vr u)_x=0,
\eeq
\beq\label{g6}
(\vr u)_t+(\vr u^2 +R\vr\vartheta+\f{1}{2}b^2)_x=\mu u_{xx},
\eeq
\beq\label{g7}
c_V(\vr\vartheta)_t+c_V(\vr u \vartheta)_x=\kappa\vartheta_{xx}-R\vr\vartheta u_x+\mu u_x^2,
\eeq
\beq\label{g8}
b_t+(b u)_x=0,
\eeq
where $\mu=2\nu+\eta$.

It should be addressed that recently Zhang and Zhao \cite{ZZ} proved the global solvability of strong solutions to the initial-boundary value problem for (\ref{g5})-(\ref{g8}) with general heat conductivity coefficient and arbitrarily large data, by making a full use of the effective viscous flux, the material derivative and the structure of the equations; while the global well-posedness of strong solutions for isentropic fluids was verified by Jiang and Zhang \cite{JZW}. When initial vacuum is allowed, Fan and Hu \cite{FH} proved the global existence of strong solutions to the problem (\ref{g5})-(\ref{g8}) under certain technical assumptions concerning the growth condition of the heat conductivity coefficient and the ratio between the  initial magnetic field and the initial density; see also \cite{YU} for a similar result dealing with the isentropic regime, where some technical hypotheses are imposed on the initial magnetic field and the initial density.

Suppose the resistivity coefficient $\lambda$ is positive. Then (\ref{g8}) reads
\beq\label{g9}
b_t+(b u)_x=\lambda b_{xx}.
\eeq
There are a great deal of results for the equations (\ref{g5})-(\ref{g7}), (\ref{g9}). In particular, Kazhikhov and Smagulov \cite{KS} announced the global well-posedness of strong solutions to the initial-boundary value problem for this system formulated in Lagrangian coordinates. Fan et al. \cite{FJGN} obtained the existence, uniqueness and Lipschitz continuous dependence on the initial data of global weak solutions to a related model. Global well-posedness of strong solutions was established by taking the effect of radiation and self-gravitation into account, see \cite{ZX}. We refer to \cite{CW1,CW2} for more results on one-dimensional model system of planar MHD.

Now we transform the MHD equations (\ref{g5})-(\ref{g8}) formulated in Eulerian coordinates to that of Lagrangian. To this end, let
\[
y=\int_0^x\vr(\xi,t)d\xi,\,\,s=t,
\]
then (\ref{g5})-(\ref{g8}) reduce to
\beq\label{f1}
\tau_t=u_x,
\eeq
\beq\label{f2}
u_t=\psi_x,
\eeq
\beq\label{f3}
c_{V}\vartheta_t=\sigma u_x+\kappa(\varrho \vartheta_x)_x,
\eeq
\beq\label{f4}
(b\tau)_t=0,
\eeq
with
\[
\tau:=\varrho^{-1}
\]
the specific volume of the fluid and
\beq\label{f5}
\sigma:=\mu\varrho u_x-R\varrho\vartheta,\,\,\psi:=\sigma-\f{1}{2}b^2.
\eeq
Here, for the purpose of convenience, we still use $(x,t)$ instead of $(y,s)$ to denote the spatial and temporal variables. Without loss of generality, we assume the conserved total mass on $[0,1]$ is one unit. Besides, the system is supplemented with the initial conditions and boundary data:
\beq\label{f6}
(\tau,u,b,\vartheta)(x,0)=(\tau_0,u_0,b_0,\vartheta_0)(x),\,\,x\in [0,1],
\eeq
\beq\label{f7}
(u,\vartheta_x)(0,t)=(u,\vartheta_x)(1,t)=0,\,\,t\in(0,\infty).
\eeq

The present paper is dedicated to obtaining global weak solutions to the initial-boundary value problem (\ref{f1})-(\ref{f7}). By establishing the boundedness of the density from above and below away from zero, weak solutions are constructed by approximation of global regular solutions, the existence of which was obtained by Zhang and Zhao \cite{ZZ} under the framework of Eulerian coordinates. Then we show the stability of weak solutions, i.e., the Lipschitz continuous dependence on the initial data. Finally, the uniqueness of global weak solutions follows as a corollary of stability result immediately. In particular, similar to the results for one-dimensional Navier-Stokes-Fourier system, see \cite{JZH, JZ} among others, our results show that neither vacuum nor concentration can form in finite time for weak solutions.

Before giving the main results of this paper, we introduce the notations used throughout this paper. Denote $\Omega:=(0,1)$, $\Omega_t:=\Omega\times(0,t)$. Let $p\in[1,\infty]$, $k$ be a positive integer. We denote the usual Lebesgue space $L^p(\Omega)$  by $L^p$, with its norm $||\cdot||_{L^p}$; $H^k$ denotes the usual Sobolev space $H^k(\Omega)$, with its norm $||\cdot||_{H^k}$; $L^p(0,T;X)$ is the space of all strongly measurable, $p$th-power integrable functions from $(0,T)$ to $X$, with $X$ being some Banach space and its corresponding norm $||\cdot||_{L^p(0,T;X)}$; $\|w\|_{V_p(\Om_t)}:=\|w\|_{L^{\infty}(0,t;{L^p})}+\|w_x\|_{L^p(0,t;{L^p})}$. $W^{1,p}(0,T;X)$ signifies the standard Sobolev space taking values in $X$. The Banach space $C([0,T];X)$ stands for all continuous functions from $[0,T]$ to $X$.

Our main results are the following two theorems.
\begin{Theorem}\label{t1}
Assume
\beq\label{f8}
u_0\in L^2,\,\vartheta_0\in L^1,\,b_0\in L^{\infty},
\eeq
\beq\label{f9}
0<m\leq\tau_0(x)\leq M <\infty,\,m\leq\vartheta_0(x),\text{ for a.e.}\,x\in \Om.
\eeq
Then there exists a weak solution $(\tau,u,b,\vartheta)$ to (\ref{f1})-(\ref{f7}) on $[0,1]\times[0,T]$ for any fixed $0<T<\infty$. Moreover, there exists a constant $C$, such that
\beq\label{f10}
C^{-1}\leq \tau(x,t) \leq C, \,C^{-1}\leq \vartheta(x,t),\text{ for a.e. }(x,t)\in \Om_T,
\eeq
\beq\label{f11}
\|u\|_{V_2(\Om_T)}+\|\tau_t\|_{L^2(0,T;L^2)}+\|b_t\|_{L^2(0,T;L^2)}\leq C,
\eeq
\beq\label{f12}
\|\vartheta\|_{V_1(\Om_T)}+\|\log\vartheta\|_{L^{\infty}(0,T;L^1)}+
\|(\log\vartheta)_x\|_{L^2(0,T;L^2)}\leq C,
\eeq
\beq\label{f15}
\|\vartheta\|_{L^{r_0}(0,T;L^{q_0})}+\|\vartheta_x\|_{L^{r_1}(0,T;L^{q_1})}\leq C,
\eeq
where
\[
\f{1}{2q_0}+\f{1}{r_0}>\f{1}{2},\,q_0,r_0\in [1,\infty],\,\f{1}{2q_1}+\f{1}{r_1}>1,\,q_1,r_1\in [1,2].
\]
\end{Theorem}
Hereafter, the same letter $C$ stands for a generic positive constant, depending only on the parameters $R, \mu, \kappa, m, M, T$ and the initial data. The definition of weak solutions will be introduced in the next section.

The next theorem gives the stability of weak solutions obtained in Theorem \ref{t1}.
\begin{Theorem}\label{t2}
Let $(\tau,u,b,\vartheta)$ and $(\widetilde{\tau},\widetilde{u},\widetilde{b},\widetilde{\vartheta})$ be two weak solutions on $[0,1]\times[0,T]$ corresponding to the initial data $(\tau_0,u_0,b_0,\vartheta_0)$ and $(\widetilde{\tau_0},\widetilde{u_0},\widetilde{b_0},\widetilde{\vartheta_0})$ satisfying (\ref{f8})-(\ref{f9}).
Then
\[
\left(\|\tau-\widetilde{\tau}\|_{L^{\infty}(\Omega_T)}+\|u-\widetilde{u}\|_{V_2(\Omega_T)}+
\|b-\widetilde{b}\|_{L^{\infty}(\Omega_T)}
+\|\vartheta-\widetilde{\vartheta}\|_{L^r(0,T;L^q)}\right)
\]
\beq\label{f13}
\leq C\left(\|\tau_0-\widetilde{\tau_0}\|_{L^{\infty}}+\|b_0-\widetilde{b_0}\|_{L^{\infty}}+
\|u_0-\widetilde{u_0}\|_{L^2}+\|\vartheta_0-\widetilde{\vartheta_0}\|_{L^1}\right),
\eeq
with arbitrary $q$ and $r$ satisfying the constraint
\beq\label{f14}
q\in[2,\infty],\,r\in[1,\infty],\,\f{1}{2q}+\f{1}{r}=\f{1}{2}\left(1+\epsilon\right),\,\,\epsilon\in\left(0,\f{1}{2}\right).
\eeq
\end{Theorem}

The uniqueness of weak solutions follows immediately from Theorem \ref{t2}.
\begin{Corollary}\label{cor2}
Let the hypotheses of Theorem \ref{t1} be satisfied. Then there exists a unique global weak solution to the initial-boundary value problem (\ref{f1})-(\ref{f7}).
\end{Corollary}
\begin{Remark}\label{rem1}
The results presented in Theorem \ref{t1}, Theorem \ref{t2} and Corollary \ref{cor2} extend the corresponding ones for the one-dimensional compressible viscous isentropic MHD equations without resistivity obtained in \cite{LYS} to the heat-conductive case.
\end{Remark}

The rest of this paper is organized as follows. In Section \ref{sec2} we recall the existence of global strong solutions due to Zhang and Zhao \cite{ZZ} under the framework of Lagrangian coordinates. In Section \ref{sec3} we prove Theorem \ref{t1} by mollifying the initial data and approximation of strong solutions. The proof of Theorem \ref{t2} is finished in Section \ref{sec4} by modifying the ideas used in \cite{AZ2,JZ,ZA1}.

\section{Global strong solutions}\label{sec2}
As pointed out in the introduction, the existence of global regular solutions to the initial-boundary value problem for (\ref{g5})-(\ref{g8}) has been proved by Zhang and Zhao \cite{ZZ}, by virtue of the effective viscous flux, the material derivative and the structure of the equations. We remark that the corresponding existence result to the initial-boundary value problem (\ref{f1})-(\ref{f7}) still holds. More precisely, the following proposition is valid.
\begin{Proposition}\label{zj}
Assume that the initial data $(\tau_0,u_0,b_0,\vartheta_0)$ satisfy
\[
(\widetilde{\tau_0},\widetilde{u_0},\widetilde{b_0},\widetilde{\vartheta_0})\in H^1,\,\, 0<m\leq \tau_0(x),\vartheta_0(x)\leq M<\infty,\text{ for any }x \in [0,1],
\]
and the compatibility condition
\[
u_0(0)=u_0(1)=0.
\]
Then there exists a unique global strong solution to the initial-boundary value problem (\ref{f1})-(\ref{f7}). Moreover, for any fixed $0<T<\infty$, valid are the estimates:
\[
C^{-1}\leq \tau(x,t),\,\vartheta(x,t) \leq C , \text{ for any }(x,t)\in[0,1]\times[0,T],
\]
\[
(\tau,b)\in L^\infty(0,T;H^1),\,\, (\tau_t,b_t)\in L^\infty(0,T;L^2),
\]
\[
(u,\vartheta)\in L^\infty(0,T;H^1)\cap L^2(0,T;H^2), \, \,(u_t,\vartheta_t)\in L^2(0,T;L^2).
\]
\end{Proposition}

The local existence of strong solution can be obtained in terms of the classical Banach fixed point theorem. Thus, it suffices to derive the global a priori estimates in order to show the global existence of strong solution. Here, for brevity, we give the following two fundamental lemmas.

First of all, the conservations of volume, total energy and entropy imply the following estimates, the proof of which is quite standard.
\begin{Lemma}\label{l1}
Let $(\tau,u,b,\vartheta)$ be a smooth solution to the initial-boundary value problem (\ref{f1})-(\ref{f7}) on $[0,1]\times[0,T]$. Then
\beq\label{a1}
\int^1_0\tau(x,t)\dx=\int^1_0\tau_0(x)\dx=1, \, \,\text{for any } t\in [0,T],
\eeq
\beq\label{a2}
\sup_{0\leq t \leq T}\int^1_0\left(\f{1}{2}u^2+c_{V}\vartheta+
\f{1}{2}b^2_0\tau^2_0\tau^{-1}\right)(x,t)\dx\leq C,
\eeq
\[
\sup_{0\leq t\leq T}\int_0^1 \left(R(\tau-\log\tau-1)+c_{V}(\vartheta-\log\vartheta-1)\right)(x,t)\dx
\]
\beq\label{a3}
+\int_0^T\int_0^1\left(\kappa \f{\vartheta_x^2}{\tau \vartheta^2}+\mu \f{u_x^2}{\tau\vartheta}\right)\dxdt \leq C.
\eeq
\end{Lemma}

Next, based on Lemma \ref{l1}, we bound the density from above and below away from zero, which plays a crucial role in the subsequent analysis. Concerning the proof, we follow and modify an idea from Antontsev et al. \cite{AKS}.
\begin{Lemma}\label{l2}
Let $(\tau,u,b,\vartheta)$ be a smooth solution to the initial-boundary value problem (\ref{f1})-(\ref{f7}) on $[0,1]\times[0,T]$. Then
\beq\label{a4}
C^{-1}\leq \varrho(x,t)\leq C,\text{ for any }(x,t)\in[0,1]\times[0,T].
\eeq
\end{Lemma}
{\bf Proof.} Using (\ref{f1}) and (\ref{f4}), we rewrite (\ref{f2}) as
\beq\label{a5}
u_t+\mu \,(\log\varrho)_{xt}+\left(R\varrho\vartheta+\f{1}{2}\f{b_0^2}{\varrho_0^2}\,\varrho^2\right)_x=0.
\eeq
Note that (\ref{a1}) and the mean value theorem show that for each $t\in [0,T]$, there exists $a(t)\in [0,1]$, such that $\varrho(a(t),t)=1$. We first integrate (\ref{a5}) over $(0,t)$ with respect to $t$, then over $(a(t),x)$ ($x$ is an arbitrary fixed point in $[0,1]$) with respect to $x$, and then take exponential on both sides of the resulting equation, finally to obtain a representation of the density:
\beq\label{a6}
\varrho(x,t)\exp\left[\f{1}{\mu}\int_0^t\left(R\varrho\vartheta+\f{1}{2}\f{b_0^2}{\varrho_0^2}\,\varrho^2\right)(x,s)\ds\right]
=\varrho_0(x)\,Y(t)\,B(x,t),
\eeq
where
\[
B(x,t):=\exp\left(\f{1}{\mu}\int_{a(t)}^x u_0(\xi)-u(\xi,t)\,d\xi\right),
\]
\[
Y(t):=\varrho_0^{-1}(a(t))\exp\left[\f{1}{\mu}\int_0^t\left(R\varrho\vartheta+\f{1}{2}\f{b_0^2}{\varrho_0^2}\,\varrho^2\right)(a(t),s)\ds\right]
.
\]
It follows easily from (\ref{a2}) and Cauchy-Schwarz's inequality that
\beq\label{a7}
C^{-1}\leq B(x,t) \leq C , \,C^{-1}\leq Y(t),\text{ for any }(x,t)\in[0,1]\times[0,T].
\eeq
Next, we need to bound $Y(t)$ from above. To this end, direct computation yields
\[
\f{\partial}{\partial t}\exp\left[\f{1}{\mu}\int_0^t\left(R\varrho\vartheta+\f{1}{2}\f{b_0^2}{\varrho_0^2}\,\varrho^2\right)(x,s)\ds\right]
\]
\[
  =  \f{1}{\mu}\left(R\varrho\vartheta+\f{1}{2}\f{b_0^2}{\varrho_0^2}\,\varrho^2\right)\exp\left[\f{1}{\mu}\int_0^t
       \left(R\varrho\vartheta+\f{1}{2}\f{b_0^2}{\varrho_0^2}\,\varrho^2\right)(x,s)\ds\right]
\]
\[
 = \f{1}{\mu}\left(R\varrho\vartheta+\f{1}{2}\f{b_0^2}{\varrho_0^2}\,\varrho^2\right)
       \varrho_0(x)\,Y(t)\,B(x,t)\varrho^{-1}(x,t)
\]
\[
 = \f{1}{\mu}\varrho_0(x)\,Y(t)\,B(x,t)\,\left(R\vartheta+\f{1}{2}\f{b_0^2}{\varrho_0^2}\,\varrho\right).
\]
Integrating the above equation over $(0,t)$ with respect to $t$ gives
\[
\exp\left[\f{1}{\mu}\int_0^t\left(R\varrho\vartheta+\f{1}{2}\f{b_0^2}{\varrho_0^2}\,\varrho^2\right)(x,s)\ds\right]
\]

\[
=1+\f{1}{\mu}\varrho_0(x)\int_0^tY(s)B(x,s)\left(R\vartheta+\f{1}{2}\f{b_0^2}{\varrho_0^2}\,\varrho\right)(x,s)\,\ds.
\]
Thus, (\ref{a6}) is reformulated as
\beq\label{a8}
\varrho(x,t)=\varrho_0(x)\,Y(t)\,B(x,t)\left[1+\f{1}{\mu}\varrho_0(x)\int_0^tY(s)B(x,s)\left(R\vartheta+\f{1}{2}\f{b_0^2}{\varrho_0^2}
\varrho\right)(x,s)\ds\right]^{-1}.
\eeq
For the convenience of later use, we write (\ref{a8}) as
\[
Y(t)\varrho^{-1}(x,t)=\left[1+\f{1}{\mu}\varrho_0(x)\int_0^tY(s)B(x,s)\left(R\vartheta+\f{1}{2}\f{b_0^2}{\varrho_0^2}\varrho\right)(x,s)\ds\right]
\varrho_0^{-1}(x).
\]
Integrating the above identity over $(0,1)$ with respect to $x$, by means of (\ref{a1}), (\ref{a2}), (\ref{a7}), yields
\beq\label{a9}
Y(t)\leq C\left(1+\int_0^tY(s)\ds\right).
\eeq
An application of Gronwall's inequality to (\ref{a9}) implies
\beq\label{a10}
Y(t)\leq C,\text{ for any } t\in[0,T].
\eeq
Consequently, we conclude from (\ref{a7}), (\ref{a8}), and (\ref{a10}) that
\beq\label{a11}
0< \varrho(x,t)\leq C.
\eeq

It remains to bound the density from below away from zero. By invoking (\ref{a7}), (\ref{a8}), (\ref{a10}) and (\ref{a11}),  we have
\beq\label{a12}
\left\|\varrho^{-1}\right\|_{L^{\infty}}\leq C\left(1+\int_0^t \|\vartheta\|_{L^{\infty}} \ds\right).
\eeq
Note that Sobolev's embedding inequality and Cauchy-Schwarz's inequality give
\[
\|\vartheta\|_{L^{\infty}}  \leq  \int_0^1\vartheta\dx+\int_0^1|\vartheta_x|\dx
\]
\[
\leq  C+\left(\int_0^1\varrho \f{\vartheta_x^2}{\vartheta^2}\dx\right)^{\f{1}{2}}\left(\int_0^1\f{\vartheta^2}{\varrho}\dx
 \right)^{\f{1}{2}}
\]
\[
\leq  C+ \left(\int_0^1\varrho \f{\vartheta_x^2}{\vartheta^2}\dx\right)^{\f{1}{2}}\left(\int_0^1\vartheta\dx\right)^{\f{1}{2}}
\left\|\varrho ^{-1}\right\|_{L^{\infty}}^{\f{1}{2}}\|\vartheta\|_{L^{\infty}}^{\f{1}{2}}
\]
\[
\leq  C+C\left(\int_0^1\varrho \f{\vartheta_x^2}{\vartheta^2}\dx\right)^{\f{1}{2}}
\left\|\varrho ^{-1}\right\|_{L^{\infty}}^{\f{1}{2}}\|\vartheta\|_{L^{\infty}}^{\f{1}{2}},
\]
where (\ref{a2}) has been used.
With the help of Cauchy-Schwarz's inequality again, we deduce from the above inequality that
\beq\label{a13}
\|\vartheta\|_{L^{\infty}}\leq C+C\,\int_0^1\varrho \f{\vartheta_x^2}{\vartheta^2}\dx \,\left\|\varrho^{-1}\right\|_{L^{\infty}}.
\eeq
Substituting (\ref{a13}) into (\ref{a12}) yields
\[
\left\|\varrho^{-1}\right\|_{L^{\infty}} \leq C+C\int_0^t \left(\int_0^1\varrho \f{\vartheta_x^2}{\vartheta^2}\dx\right)\,
\left\|\varrho^{-1}\right\|_{L^{\infty}} \ds,
\]
which, by recalling (\ref{a3}) and using Gronwall's inequality, gives
\[
\left\|\varrho^{-1}\right\|_{L^{\infty}} \leq C,\text{ for any }t\in [0,T].
\]
This completes the proof of Lemma \ref{l2}.

\begin{Remark}\label{rem1}
Integrating (\ref{a13}) over $(0,T)$ with respect to $t$, on account of (\ref{a3}) and (\ref{a4}), gives the bound
\[
\int_0^T \|\vartheta\|_{L^{\infty}}\dt \leq C,
\]
which combined with (\ref{a2}) and the interpolation formula show that
\beq\label{a14}
\|\vartheta\|_{L^p(0,T;L^q)}\leq C,  \text{ for any }p,q\in [1,\infty],\,\f{1}{p}+\f{1}{q}=1.
\eeq
As a consequence, (\ref{a3}), (\ref{a4}) and (\ref{a14}) allow us to estimate
\[
\int_0^T\int_0^1 |\vartheta_x| \dxdt  \leq  \left(\int_0^T\int_0^1\varrho \, \f{\vartheta_x^2}{\vartheta^2}\,\dxdt\right)^\f{1}{2}\left(\int_0^T\int_0^1\f{\vartheta^2}
{\varrho}\,\dxdt\right)^\f{1}{2}
\]
\[
 \leq  C \,\left(\int_0^T\int_0^1\vartheta^2 \dxdt\right)^{\f{1}{2}} \,\,\leq  C ,
\]
where Cauchy-Schwarz's inequality has been invoked.
\end{Remark}
\begin{Remark}\label{rem3}
It should be noticed that in \cite{ZZ} the authors proved the boundedness of density from above and below by making a full use of the material derivative and the structure of the equations in the Eulerian coordinates. In particular, the lower boundedness of the density follows as a consequence of the boundedness of the magnetic field. In our case, the boundedness of magnetic field follows directly from Lemma \ref{l2}.
\end{Remark}

Based on Lemmas \ref{l1}-\ref{l2}, we can further derive the necessary global a priori estimates for $(\tau, u, b, \vartheta)$. Here we only list those without giving the detailed proof (see \cite{AKS,ZZ}).
\begin{Lemma}\label{l3}
Let $(\tau,u,b,\vartheta)$ be a smooth solution to the initial-boundary value problem (\ref{f1})-(\ref{f7}) on $[0,1]\times[0,T]$. Then
\[
\sup_{0\leq t\leq T}\int_0^1\left(\tau_x^2+b_x^2\right)(x,t)\dx \leq C,
\]
\[
\sup_{0\leq t\leq T}\int_0^1u_x^2(x,t) \dx + \int_0^T\int_0^1\left(u_{xx}^2+u_t^2\right)(x,t)\dxdt \leq C,
\]
\[
\sup_{0\leq t\leq T}\int_0^1\left(\tau_t^2+b_t^2\right)(x,t)\dx \leq C,
\]
\[
\sup_{0\leq t\leq T}\int_0^1\vartheta_x^2(x,t) \dx + \int_0^T\int_0^1\left(\vartheta_{xx}^2+\vartheta_t^2\right)(x,t)\dxdt \leq C,
\]
\[
C^{-1}\leq \vartheta(x,t)\leq C,\text{ for any }(x,t)\in[0,1]\times[0,T].
\]
\end{Lemma}

As a consequence, we complete the proof of Proposition \ref{zj} by extending the local solutions globally with respect to time, making use of the global a priori estimates obtained in Lemmas \ref{l1}-\ref{l3}.
Before ending this section, we introduce the definition of weak solution to the MHD system (\ref{f1})-(\ref{f4}).
\begin{Definition}\label{dfn1}
A quadruple $(\tau,u,b,\vartheta)$ is said to be a weak solution to the MHD system (\ref{f1})-(\ref{f4}) on $[0,1]\times [0,T]$ , where $0<T<\infty$ is arbitrarily fixed, with boundary conditions (\ref{f7}) and initial data $(\tau_0,u_0,b_0,\vartheta_0)$ satisfying
\[
\tau_0 ,\tau^{-1}_0\in L^{\infty},  \tau_0>0,  u_0\in L^2, b_0\in L^{\infty}, \vartheta_0 \in L^1,
\vartheta_0^{-1} \in L^{\infty},\vartheta_0 >0,
\]
provided that
\[
\tau \in W^{1,2}(0,T;L^2), \,u\in L^\infty(0,T;L^2)\cap L^2(0,T;H^1_0), \,b\in W^{1,2}(0,T;L^2),
\]
\[
\vartheta \in L^{\infty}(0,T;L^1),\,\,\vartheta_x \in L^1(0,T;L^1),
\]
\[
\tau,\tau^{-1}\in L^{\infty}(\Omega_T),\, \tau>0,\,\vartheta^{-1} \in L^{\infty}(\Omega_T),\,\vartheta >0,
\]
\[
\tau_t=u_x \, \, \text{a.e. in } \Om_T, \, \, \tau(x,0)=\tau_0(x)  \text { for a.e. } x \in \Om,
\]
\[
b=b_0 \tau_0 \tau^{-1}  \text{ a.e. in }\Om_T,
\]
and the following integral identities hold:
\beq\label{dfn2}
\int_0^T \int_{\Omega} u\phi_t-\left(\mu \varrho u_x-R\varrho\vartheta-\f{1}{2}b^2\right)\phi_x \,\dxdt+\int_{\Om}u_0\phi(x,0)\dx=0,
\eeq
\beq\label{dfn3}
\int_0^T \int_{\Omega}c_V \vartheta \varphi_t-\kappa \varrho \vartheta_x \varphi_x+\sigma u_x\varphi\,\dxdt+\int_{\Om}c_V \vartheta_0\varphi(x,0)\dx=0,
\eeq
for any $\phi \in C^{\infty}_c(\Omega\times [0,T))$ and $\varphi \in C^{\infty}_c(\bar{\Omega} \times [0,T))$.
\end{Definition}

\section{Existence of global weak solutions}\label{sec3}
This section is designed to give the proof of Theorem \ref{t1} by mollifying the initial data and approximation of strong solutions, the existence of which is guaranteed by Proposition \ref{zj}. Then the existence of weak solutions is established in terms of the analysis of weak convergence.

We begin with the construction of strong solutions. Let $(\tau_0,u_0,b_0,\vartheta_0)$ be subject to conditions (\ref{f8})-(\ref{f9}), then, by regularizing the initial data, there exists $(\tau_0^{\epsilon},u_0^{\epsilon},b_0^{\epsilon},\vartheta_0^{\epsilon})$ for each $\epsilon>0$, such that
\[
\tau_0^{\epsilon},\,b_0^{\epsilon},\,\vartheta_0^{\epsilon}\,\in C^2([0,1]),\,\,u_0^{\epsilon} \in C_c^2(\Om),
\]
\[
C^{-1}\leq \tau_0^{\epsilon}\leq C,\,\,C^{-1}\leq \vartheta_0^{\epsilon},\,\,\|b_0^{\epsilon}\|_
{L^{\infty}}\leq \|b_0\|_{L^{\infty}},
\]
\[
(\tau_0^{\epsilon},u_0^{\epsilon},b_0^{\epsilon}) \,\rightarrow(\tau_0,u_0,b_0)\text{ strongly in } L^2,\,\vartheta_0^{\epsilon}\rightarrow \vartheta_0\,\,\text{strongly in}\,L^1,\text{ as }\epsilon \rightarrow 0^{+}.
\]
It follows from Proposition \ref{zj} that there exists a unique global strong solution $(\tau^{\epsilon},u^{\epsilon},b^{\epsilon},\vartheta^{\epsilon})$ to the initial-boundary value problem (\ref{f1})-(\ref{f7}) with $(\tau_0,u_0,b_0,\vartheta_0)$ replaced by $(\tau_0^{\epsilon},u_0^{\epsilon},b_0^{\epsilon},\vartheta_0^{\epsilon})$. Moreover, Lemmas \ref{l1}-\ref{l2} and Remark \ref{rem1} give the following uniform-in-$\epsilon$ bounds:
\beq\label{b1}
C^{-1}\leq \tau^{\epsilon}(x,t) \leq C , \,\,C^{-1}\leq \vartheta^{\epsilon}(x,t)\text{ for any }(x,t)\in[0,1]\times[0,T],
\eeq
\beq\label{b2}
\|u^{\epsilon}\|_{L^\infty(0,T;L^2)}+ \|u^{\epsilon}_x\|_{L^2(0,T;L^2)}\leq C,
\eeq
\beq\label{b3}
\|\tau^{\epsilon}_t\|_{L^2(0,T;L^2)}+\|b^{\epsilon}_t\|_{L^2(0,T;L^2)}\leq C,
\eeq
\beq\label{b4}
\|\vartheta^{\epsilon}\|_{L^{\infty}(0,T;L^1)}+\|\vartheta_x^{\epsilon}\|_{L^1(0,T;L^1)}+\|\log \vartheta^{\epsilon}\|_{L^{\infty}(0,T;L^1)}+
\|(\log\vartheta^{\epsilon})_x\|_{L^2(0,T;L^2)}\leq C,
\eeq
\beq\label{b5}
\|\vartheta^{\epsilon}\|_{L^p(0,T;L^q)}\leq C,  \text{ for any }p,q\in [1,\infty],\,\f{1}{p}+\f{1}{q}=1.
\eeq
In addition, due to Lemma 8 in \cite{AZ2}, valid is the uniform bound 
\beq\label{b6}
\|\vartheta^{\epsilon}_x\|_{L^{r_1}(0,T;L^{q_1})}\leq C,  \text{ for any }q_1,r_1\in [1,2],\,\f{1}{2q_1}+\f{1}{r_1}>1.
\eeq

Using (\ref{b1})-(\ref{b6}), we can extract a subsequence of $(\tau^{\epsilon},u^{\epsilon},b^{\epsilon},\vartheta^{\epsilon})$, still denoted by $(\tau^{\epsilon},u^{\epsilon},b^{\epsilon},\vartheta^{\epsilon})$, such that as $\epsilon\rightarrow0^{+}$, there holds:
\[
\tau^{\epsilon}\rightarrow\tau \text{ weakly}-\star \text{ in }L^{\infty}(0,T;L^{\infty}),
\]
\beq\label{tb2}
u^{\epsilon}\rightarrow u  \text{ weakly}-\star \text{ in }L^{\infty}(0,T;L^2),
\eeq
\beq\label{tb3}
(\tau^{\epsilon}_t,u^{\epsilon}_x)\rightarrow(\tau_t,u_x) \text{ weakly in }L^2(0,T;L^2),
\eeq
\[
\vartheta^{\epsilon}\rightarrow \vartheta  \text{ weakly or weakly}-\star \text{ in }L^{p}(0,T;L^q),\text{ with }p,\,q\text{ satisfying }(\ref{b5}),
\]
\beq\label{tb1}
\vartheta^{\epsilon}_x\rightarrow \vartheta_x  \text{ weakly in } L^{r_1}(0,T;L^{q_1}),\text{ with }r_1,\,q_1\text{ satisfying }(\ref{b6}).
\eeq
Furthermore, by the theorem of compactness embedding, there holds:
\[
u^{\epsilon}\rightarrow u  \text{ strongly in }L^2(0,T;L^2),
\]
\[
\vartheta^{\epsilon}\rightarrow \vartheta \text{ strongly in } L^1(0,T;L^2).
\]

In view of the nonlinearity of the equations, we are now in a position to show that $\tau^{\epsilon}$ converges to $\tau$ strongly as $\epsilon\rightarrow0^{+}$. The next lemma plays an important role in the analysis of compactness.
\begin{Lemma}\label{l4}
For any $0<h<1$, there holds
\beq\label{b7}
\|\Delta_h\tau^{\epsilon}\|_{L^\infty(0,T;L^2)}\leq C(\|\Delta_h \tau_0\|_{L^2}+\|\Delta_h b_0\|_{L^2}+h).
\eeq
\end{Lemma}
{\bf Proof.} For the sake of simplicity, we set
\[
 \psi^{\epsilon}:=\mu \varrho^{\epsilon} u^{\epsilon}_x-R\varrho^{\epsilon}\vartheta^{\epsilon}-\f{1}{2}(b^{\epsilon})^2,\,\,
 a^{\epsilon}:=b_0^{\epsilon}\tau_0^{\epsilon},
\]
\[
\sigma^{\epsilon}:=\mu \varrho^{\epsilon} u^{\epsilon}_x-R\varrho^{\epsilon}\vartheta^{\epsilon},\,\,
B^{\epsilon}(x,t):=\exp\left(\f{1}{\mu}\int_0^t\psi^{\epsilon}(x,s)\ds\right),
\]
then $(\tau^{\epsilon},u^{\epsilon},b^{\epsilon},\vartheta^{\epsilon})$ satisfies the system below:
\beq\label{b8}
\tau^{\epsilon}_t=u^{\epsilon}_x,
\eeq
\beq\label{b9}
u^{\epsilon}_t= \psi^{\epsilon}_x,
\eeq
\beq\label{b10}
c_{V}\vartheta^{\epsilon}_t=\sigma^{\epsilon}u^{\epsilon}_x+\kappa(\varrho ^{\epsilon}\vartheta^{\epsilon}_x)_x,
\eeq
\beq\label{b11}
(b^{\epsilon}\tau^{\epsilon})_t=0.
\eeq
Obviously, by virtue of (\ref{b1}) and (\ref{b5}), the following estimate is valid:
\beq\label{b12}
C^{-1}\leq B^{\epsilon}(x,t) \leq C, \text{ for any }(x,t)\in[0,1]\times[0,T].
\eeq
Note that (\ref{b8}), (\ref{b9}) and (\ref{b11}) together imply
\beq\label{b13}
\tau^{\epsilon}_t=\f{\psi^{\epsilon}}{\mu}\tau^{\epsilon}+\f{R}{\mu}\vartheta^{\epsilon}+
\f{1}{2\mu}(a^{\epsilon})^2 (\tau^{\epsilon})^{-1}
\eeq
Multiplying (\ref{b13}) by $\exp\left(-\f{1}{\mu}\int_0^t\psi^{\epsilon}(x,s)\ds\right)$ and integrating the resulting equation over $(0,t)$ with respect to $t$ gives
\[
\tau^{\epsilon}=B^{\epsilon}\left(\tau^{\epsilon}_0+\f{1}{\mu}
\int_0^t \f{(R\vartheta^{\epsilon}+\f{1}{2}( a^{\epsilon})^2(\tau^{\epsilon})^{-1})(x,\xi) }{B^{\epsilon}(x,\xi)}d\xi \right),
\]
and one easily calculates that
\[
\Delta_h\tau^{\epsilon}(x,t)= \Delta_h B^{\epsilon}(x,t)\left[\tau^{\epsilon}_0(x+h)+\f{1}{\mu}
                             \int_0^t \f{(R\vartheta^{\epsilon}+\f{1}{2}( a^{\epsilon})^2(\tau^{\epsilon})^{-1})(x+h,\xi) }{B^{\epsilon}(x+h,\xi)}d\xi\right]
\]
\[
                             +  B^{\epsilon}(x,t)\left(\Delta_h\tau^{\epsilon}_0+\f{1}{\mu}\int_0^t \f{\Delta_h (R\vartheta^{\epsilon}+\f{1}{2}( a^{\epsilon})^2(\tau^{\epsilon})^{-1})(x,\xi) }{B^{\epsilon}(x+h,\xi)}d\xi\right)
\]
\[
                             -  B^{\epsilon}(x,t)\left( \f{1}{\mu} \int_0^t \f{(R\vartheta^{\epsilon}+\f{1}{2}( a^{\epsilon})^2(\tau^{\epsilon})^{-1})(x,\xi) }{B^{\epsilon}(x,\xi)B^{\epsilon}(x+h,\xi)}\Delta_h B^{\epsilon}(x,\xi) d\xi\right),
\]
By virtue of (\ref{b1}), (\ref{b5}), (\ref{b6}) and (\ref{b12}), we estimate
\[
\|\Delta_h\tau^{\epsilon}(x,t)\|_{L^2} \leq  C\left(\|\Delta_hB^{\epsilon}(x,t)\|_{L^\infty(0,T;L^2)}+
                                       \|\Delta_h\tau^{\epsilon}_0\|_{L^2}\right)
\]
\[
                                       + C \int_0^t\|\Delta_h \vartheta^{\epsilon}(x,\xi)\|_{L^2}
                                       +\|\Delta_h\tau^{\epsilon}(x,\xi)\|_{L^2}
                                       +\|\Delta_h a^{\epsilon}\|_{L^2} d\xi
\]
\[
 \leq  C\left(h\|u^{\epsilon}-u_0^{\epsilon}\|_{L^\infty(0,T;L^2)}+h\|\vartheta^{\epsilon}_x\|_{L^1(0,T;L^2)}\right)
 \]
 \[
 +C\left(\|\Delta_h\tau^{\epsilon}_0\|_{L^2}
  + \|\Delta_h b^{\epsilon}_0\|_{L^2}
  +\int_0^t\|\Delta_h\tau^{\epsilon}(x,\xi)\|_{L^2} d\xi\right)
\]
\[
\leq C\left(h+\|\Delta_h\tau_0\|_{L^2}
  + \|\Delta_h b_0\|_{L^2}+\int_0^t\|\Delta_h\tau^{\epsilon}(x,\xi)\|_{L^2} d\xi\right),
\]
then (\ref{b7}) follows from the above inequality immediately with an application of Gronwall's inequality. This completes the proof of Lemma \ref{l4}.

Combining  (\ref{b3}) and Lemma \ref{l4}, we obtain that for any $0<h<1$, $0<s<T$,
\beq\label{b14}
\|\tau^{\epsilon}(\cdot+h,\cdot+s)-\tau^{\epsilon}\|_{L^{\infty}(0,T-s;L^2)}
\leq C\left(\|\Delta_h \tau_0\|_{L^2}+\|\Delta_h b_0\|_{L^2}+h+s^{\f{1}{2}}\right).
\eeq
Making use of the criterion of compactness of sets in $L^2(0,T;L^2)$, (\ref{b14}) gives
\beq\label{b15}
\tau^{\epsilon}\rightarrow \tau\text{ strongly in } \,L^2(0,T;L^2)\text{ as }\epsilon\rightarrow0^{+},
\eeq
and furthermore, by invoking (\ref{b1}),
\beq\label{b16}
\tau^{\epsilon}\rightarrow \tau \text{ strongly in }L^p(0,T;L^p)\text{ as }\epsilon\rightarrow0^{+},\text{ for any }1 \leq p <\infty.
\eeq
Defining
\[
b:=b_0\tau_0\tau^{-1},
\]
then one deduces from (\ref{b1}) and (\ref{b15}) that
\beq\label{b17}
b^{\epsilon}\rightarrow b \text{ strongly in }L^2(0,T;L^2)\text{ as }\epsilon\rightarrow0^{+},
\eeq
\beq\label{b18}
b^{\epsilon}_t\rightarrow b_t\text{ weakly in } L^2(0,T;L^2)\text{ as }\epsilon\rightarrow0^{+},
\eeq
\beq\label{b19}
\|b_t\|_{L^2(0,T;L^2)}\leq C.
\eeq

The proof of Theorem \ref{t1} follows essentially the same line as that of \cite{AZ2,JZ}. Specifically, estimates (\ref{f10})-(\ref{f12}) are easily obtained from the uniform-in-$\epsilon$ bounds (\ref{b1})-(\ref{b6}) with an application of the theorem of Fatou. Thus, it remains to show that integral identities (\ref{dfn2})-(\ref{dfn3}) hold. Note that (\ref{b1}), (\ref{b2}) and (\ref{b5}) give the bound
\[
\|\psi^{\epsilon}\|_{L^2(0,T;L^2)} \leq C,
\]
which, by using (\ref{b1}), (\ref{b2}), (\ref{b5}), (\ref{b16}), yields
\beq\label{b20}
\psi^{\epsilon}\rightarrow \psi\text{ weakly in } L^2(0,T;L^2)\text{ as }\epsilon\rightarrow0^{+}.
\eeq
Therefore, we multiply (\ref{b9}) by any $\phi \in C^{\infty}_c(\Omega\times [0,T))$ and integrate by parts, letting $\epsilon\rightarrow 0^+$, to find (\ref{dfn2}) follows. Next, we have (see \cite{AZ2})
\beq\label{b21}
\sigma^{\epsilon}u^{\epsilon}\rightarrow \sigma u \text{ weakly in } L^{\f{3}{2}-\delta}(\Om_T) \text{ for any } \delta \in (0,\f{1}{2}).
\eeq
By means of (\ref{tb2}), (\ref{tb3}), (\ref{tb1}), (\ref{b16}), (\ref{b17}), (\ref{b21}), we first multiply (\ref{b9}) by $u$ and add the resulting equation to (\ref{b10}), followed by testing the equation with any $\varphi \in C^{\infty}_c(\bar{\Omega} \times [0,T))$, integrating by parts and letting $\epsilon\rightarrow 0^+$, to obtain
\[
\int_0^T\int_{\Om} (\f{1}{2}u^2+c_V\vartheta)\varphi_t-\sigma u \varphi_x -\kappa\varrho\vartheta_x\varphi_x+\f{1}{2}b^2(u_x\varphi+u\varphi_x)\dxdt
\]
\beq\label{b22}
+\int_{\Om}(\f{1}{2}u_0^2+c_V\vartheta_0)
\varphi(x,0)\dx=0.
\eeq
To show the validity of (\ref{dfn3}), we approximate the limit functions $(u,\psi)$ by $(u_{\epsilon},\psi_{\epsilon})$, with a suitable regularization. Moreover, we deduce from (\ref{dfn2}), upon an appropriate choice of $\phi$,
\[
(u_{\epsilon})_t=(\psi_{\epsilon})_x \text{ in }L^2(\Om_T).
\]
By testing the above equality with $u_{\epsilon}\varphi$, we conclude that, with an application of integration by parts and letting $\epsilon\rightarrow 0^+$,
\[
\int_0^T\int_{\Om}-\f{1}{2}u^2\varphi_t+\psi(u_x\varphi+u\varphi_x)\dxdt=\int_{\Om}\f{1}{2}u_0^2\varphi(x,0)\dx.
\]
Thus (\ref{dfn3}) follows immediately by subtracting the above integral identity from (\ref{b22}). The proof of Theorem \ref{t1} is finished.

\section{Stability of weak solutions}\label{sec4}
In this section, we give the proof of Theorem \ref{t2} by virtue of modifying the ideas adopted in  \cite{AZ2,JZ,ZA1}. To begin with, we have the following lemma, the proof of which is omitted.
\begin{Lemma}\label{lem1}
Let the assumptions of Theorem \ref{t2} be satisfied. Then the following representations are valid in $(0,1)\times (0,T)$:
\[
\tau(x,t) = \exp\left(\f{1}{\mu}\int_0^t\psi(x,s)\ds\right)
\]
\beq\label{kb1}
          \times  \left[ \tau_0+\int_0^t\exp\left(-\f{1}{\mu}\int_0^{\xi}\psi(x,s)\ds\right)
          \left(\f{R}{\mu}\vartheta+\f{1}{2\mu}b_0^2\tau_0^2\tau^{-1}\right)(x,\xi) d\xi \right],
\eeq
and
\beq\label{kb2}
\int_0^t\psi(x,s)\ds=(J_{\Omega}(u-u_0))(x,t)+\int_0^t<\psi(\cdot,s)> \ds,
\eeq
where the linear operator $J_{\Omega}$ is defined by
\[
J_{\Omega}w(x):=\int_0^x w(\xi)d\xi-<\int_0^x w(\xi)d\xi>,\,\, <w>:=\int_0^1w(x)\dx.
\]
\end{Lemma}

In order to simplify the expressions below, we introduce the abbreviated notations.
\[
(\Delta \tau,\Delta u, \Delta b,\Delta \vartheta):=(\tau-\widetilde{\tau},u-\widetilde{u},b-\widetilde{b},\vartheta-\widetilde{\vartheta}),
\]
\[
(\Delta \tau_0,\Delta u_0, \Delta b_0,\Delta \vartheta_0):=(\tau_0-\widetilde{\tau_0},u_0-\widetilde{u_0},b_0-\widetilde{b_0}, \vartheta_0-\widetilde{\vartheta_0}),
\]
\[
\Delta\sigma:=\sigma-\widetilde{\sigma},\,
\widetilde{\sigma}:=\mu \widetilde{\varrho}\widetilde{u}_x-R\widetilde{\varrho}\widetilde{\vartheta},
\]
\[
\Delta \psi=\psi-\widetilde{\psi},\,\widetilde{\psi}:=\widetilde{\sigma}-
\f{1}{2}\widetilde{b}^2,
\]
\[
g:=\exp\left(\f{1}{\mu}\int_0^t\psi(x,s)\ds\right),\,\,
\widetilde{g}:=\exp\left(\f{1}{\mu}\int_0^t\widetilde{\psi}(x,s)\ds\right),
\]
\[
K:=\f{R}{\mu}\vartheta+\f{1}{2\mu}b_0^2\tau_0^2\tau^{-1},\,\,
\widetilde{K}:=\f{R}{\mu}\widetilde{\vartheta}+\f{1}{2\mu}(\widetilde{b_0})^2(\widetilde{\tau_0})^2
(\widetilde{\tau})^{-1},
\]
\[
\widetilde{\vr}:=(\widetilde{\tau})^{-1},\,\,\Delta\vr:=\vr-\widetilde{\vr}.
\]

The next lemma concerning the supremum of $\Delta\tau$ plays a crucial role for the proof of Theorem \ref{t2}.
\begin{Lemma}\label{lem2}
Let the assumptions of Theorem \ref{t2} be satisfied. Then
\[
\|\Delta\tau\|_{L^{\infty}(\Omega_t)} \leq C(\|\Delta \tau_0\|_{L^{\infty}}+\|\Delta b_0\|
                                      _{L^{\infty}}+\|\Delta u_0\|_{L^2}
\]
\beq\label{kb3}
+\|\Delta u\|_{V_2(\Om_t)}
+\|\Delta\vartheta\|_{L^1(0,t;L^{\infty})},\text{ for any } t\in(0,T].
\eeq
\end{Lemma}
{\bf Proof.} In accordance with the representation (\ref{kb1}), we see
\beq\label{kb4}
\Delta\tau  = g\left[\Delta \tau_0+\int_0^t K\left(\f{1}{g}-\f{1}{\widetilde{g}}\right)+\f{K-\widetilde{K}}
           {\widetilde{g}}d\xi\right]
           + (g-\widetilde{g})\left(\widetilde{\tau_0}+\int_0^t\f{\widetilde{K}}{\widetilde{g}}d\xi\right).
\eeq
On account of the regularity results stated in Theorem \ref{t2}, we deduce that
\[
\psi \in L^2(\Om_T),\,\left(\int_0^t\psi(x,s)\ds\right)_x \in L^{\infty}(0,T;L^2),
\]
and moreover
\[
\int_0^t\psi(x,s)\ds \in C(\bar{\Om}_T).
\]
Consequently,
\beq\label{kb5}
C^{-1}\leq g,\,\,\widetilde{g} \leq C.
\eeq
Using (\ref{f10}) and (\ref{kb5}), we can estimate (\ref{kb4}) through
\[
|\Delta\tau|\leq C\left(|\Delta \tau_0|+\int_0^t(1+\vartheta)\left|\int_0^{\xi}\Delta\psi\ds\right|+|\Delta\tau|+|\Delta b_0|+|\Delta \tau_0| +|\Delta\vartheta| d\xi\right)
\]
\[
+C\left|\int_0^t\Delta\psi\ds\right|\left(1+\int_0^t\widetilde{\vartheta} d\xi\right),
\]
which by virtue of (\ref{f15}) shows
\[
\|\Delta\tau(\cdot,t)\|_{L^{\infty}}\leq C\left(\|\Delta \tau_0\|_{L^{\infty}}
                                    +\|\Delta b_0\|_{L^{\infty}}+\left\|\int_0^{\xi}\Delta\psi \ds\right\|_{L^{\infty}
                                    (\Omega_t)}
                                    +\int_0^t\|\Delta\tau(\cdot,\xi)\|_{L^{\infty}} d\xi \right)
\]
\beq\label{kb6}
+C\|\Delta\vartheta\|_{L^1(0,t;L^{\infty})}.
\eeq
By means of (\ref{kb2}), we have
\beq\label{kb7}
\left\|\int_0^{\xi}\Delta\psi \ds\right\|_{L^{\infty}(\Omega_t)}  \leq  \|J_{\Omega}\Delta u_0\|_{L^{\infty}}
                                                        +\|J_{\Omega}\Delta u\|_
                                                        {L^{\infty}(\Omega_t)}
                                                        + \|\Delta\psi\|_{L^1(\Omega_t)}.
\eeq
Obviously, there holds
\beq\label{kb8}
\Delta\psi=\mu\tau^{-1}(\Delta u)_x+\chi,
\eeq
where
\[
\chi:=\mu(\Delta \vr)\widetilde{u}_x-\mu\f{K}{\tau}+\mu\f{\widetilde{K}}{\widetilde{\tau}}.
\]
Moreover, by using (\ref{f10}), we obtain
\[
|\Delta\psi|\leq C |(\Delta u)_x|+|\chi|,
\]
\beq\label{kb9}
|\chi|\leq C |\Delta\tau|(|\widetilde{u}_x|+\widetilde{\vartheta}+1)+C (|\Delta b_0|+|\Delta \tau_0|+|\Delta \vartheta|).
\eeq
As a consequence,
\[
\|\Delta\psi\|_{L^1(\Omega_t)}\leq C\left(\|(\Delta u)_x\|_{L^1(\Omega_t)}+\|\Delta b_0\|_{L^{\infty}}+\|\Delta \tau_0\|_{L^{\infty}}+\|\Delta \vartheta\|_{L^1(\Om_t)}\right)
\]
\beq\label{kb10}
+\int_0^t \zeta(s)\|\Delta\tau(\cdot,s)\|_{L^{\infty}}\ds,
\eeq
where
\[
\zeta(t):=\|\widetilde{u}_x(\cdot,t)\|_{L^2}+\|\widetilde{\vartheta}(\cdot,t)\|_{L^2}+1.
\]
In view of (\ref{f11}) and (\ref{f15}), we see
\[
\|\zeta\|_{L^2(0,T)}\leq C.
\]
By taking advantage of (\ref{kb7}) and (\ref{kb10}), we  furthermore estimate (\ref{kb6}) as follows:
\[
\|\Delta\tau(\cdot,t)\|_{L^{\infty}}\leq C(\|\Delta \tau_0\|_{L^{\infty}}
                                    +\|\Delta b_0\|_{L^{\infty}}+\|J_{\Omega}\Delta u_0\|_{L^{\infty}}
                                    +\|J_{\Omega}\Delta u\|_ {L^{\infty}(\Omega_t)})
\]
\beq\label{kb11}
  +C\left(\|(\Delta u)_x\|_{L^1(\Omega_t)}+\|\Delta \vartheta\|_{L^1(0,t;L^{\infty})}+
                                    \int_0^t \overline{\zeta}(s)\|\Delta\tau(\cdot,s)\|_{L^{\infty}}\ds\right),
\eeq
where we have set
\[
\overline{\zeta}(t):=\zeta(t)+1.
\]
Thus (\ref{kb3}) follows immediately from (\ref{kb11}) by applying Gronwall's inequality and recalling the definition of $J_\Om$.

We now give a lemma concerning the energy estimate of $\Delta u$.
\begin{Lemma}\label{lem3}
Let the hypotheses of Theorem \ref{t2} be satisfied. Then for any $t\in (0,T]$,
\[
\|\Delta u\|_{V_2(\Om_t)} \leq C(\|\Delta u_0\|_{L^2} + \|\Delta b_0\|_{L^{\infty}}+ \|\Delta \tau_0\|_{L^{\infty}}
\]
\beq\label{kb12}
+\|\Delta \vartheta\|_{L^2(\Om_t)}+ \|\zeta\|\Delta\tau
(\cdot,s)\|_{L^{\infty}} \|_{L^2(0,t)}).
\eeq
Here $\zeta(t)=\|\widetilde{u}_x(\cdot,t)\|_{L^2}+\|\widetilde{\vartheta}(\cdot,t)\|_{L^2}+1$.
\end{Lemma}
{\bf Proof.}
By invoking (\ref{kb8}) and (\ref{f2}), we see
\[
(\Delta u)_t=[\mu \tau^{-1}(\Delta u)_x+\chi]_x.
\]
Standard energy estimate for linear parabolic equation shows that
\beq\label{kb13}
\|\Delta u\|_{V_2(\Om_t)} \leq C(\|\Delta u_0\|_{L^2} +\|\chi\|_{L^2(\Omega_t)}),
\eeq
where (\ref{f10}) has been used. On the other hand, we deduce from (\ref{kb9}) that
\beq\label{kb14}
\|\chi\|_{L^2(\Omega_t)}\leq C(\|\zeta\|\Delta\tau(\cdot,s)\|_{L^{\infty}} \|_{L^2(0,t)}+
\|\Delta b_0\|_{L^{\infty}}+ \|\Delta \tau_0\|_{L^{\infty}}+\|\Delta \vartheta\|_{L^2(\Om_t)}).
\eeq
Thus we finish the proof by combining (\ref{kb13}) and (\ref{kb14}).

At this stage, it is a routine matter to rewrite the estimates obtained in Lemmas \ref{lem2}-\ref{lem3} together. To be more precise, with an application of Gronwall's inequality, there holds
\begin{Lemma}\label{lem4}
Let the hypotheses of Theorem \ref{t2} be satisfied. Then for any $t\in (0,T]$,
\[
\|\Delta\tau\|_{L^{\infty}(\Omega_t)}+\|\Delta u\|_{V_2(\Om_t)} \leq C(\|\Delta \tau_0\|_{L^{\infty}}+\|\Delta b_0\|_{L^{\infty}}+\|\Delta u_0\|_{L^2})
\]
\beq\label{kb15}
+C(\|\Delta\vartheta\|_{L^1(0,t;L^{\infty})}+\|\Delta \vartheta\|_{L^2(\Om_t)}).
\eeq
\end{Lemma}

To proceed, we report a key lemma that concerns an estimate for $\Delta\vartheta$. We refer to Lemma 3.5 in \cite{ZA1} for the detailed proof and \cite{JZ} dealing with the case of Cauchy problem.
\begin{Lemma}\label{lem5}
Let the hypotheses of Theorem \ref{t2} be satisfied. Then for any $t\in (0,T]$,
\beq\label{kb16}
\|\Delta{\vartheta}\|_{L^r(0,t;L^q)}\leq C(\|\Delta\tau\|_{L^{\infty}(\Omega_t)}+\|\Delta u\|_{V_2(\Om_t)}+\|\Delta \vartheta\|_{L^2(\Om_t)}+\|\Delta \tau_0\|_{L^{\infty}}+\|\Delta \vartheta_0\|_{L^1}),
\eeq
where $q$ and $r$ are indicated in Theorem \ref{t2}.
\end{Lemma}

With Lemmas \ref{lem1}-\ref{lem5} at hand, we now give the proof of Theorem \ref{t2}.

{\bf Proof of Theorem \ref{t2}.} Multiplying (\ref{kb16}) by $\f{1}{2C}$ and adding the resulting inequality to (\ref{kb15}), we find
\beq\label{kb17}
\|\Delta\tau\|_{L^{\infty}(\Omega_t)}+\|\Delta u\|_{V_2(\Om_t)}+\|\Delta{\vartheta}\|_{L^r(0,t;L^q)}\leq C (c_0+\|\Delta\vartheta\|_{L^1(0,t;L^{\infty})}+\|\Delta \vartheta\|_{L^2(\Om_t)}),
\eeq
where we have set
\[
c_0:=\|\Delta \tau_0\|_{L^{\infty}}+\|\Delta b_0\|_{L^{\infty}}+\|\Delta u_0\|_{L^2}+\|\Delta \vartheta_0\|_{L^1}.
\]
Moreover, it follows from (\ref{f10}) and Definition \ref{dfn1} that
\beq\label{kb18}
\|\Delta b\|_{L^{\infty}(\Omega_t)}\leq C(\|\Delta \tau_0\|_{L^{\infty}}+\|\Delta b_0\|_{L^{\infty}}
+\|\Delta\tau\|_{L^{\infty}(\Omega_t)}),
\eeq
Therefore, we multiply (\ref{kb18}) by $\f{1}{2C}$ and add the resulting equation to (\ref{kb17}) to deduce that
\[
\|\Delta\tau\|_{L^{\infty}(\Omega_t)}+\|\Delta u\|_{V_2(\Om_t)}+\|\Delta{\vartheta}\|_{L^r(0,t;L^q)}+\|\Delta b\|_{L^{\infty}(\Omega_t)}
\]
\beq\label{kb19}
\leq C (c_0+\|\Delta\vartheta\|_{L^1(0,t;L^{\infty})}+\|\Delta \vartheta\|_{L^2(\Om_t)}),
\eeq
In addition, by virtue of choosing $q=\infty,2$, we conclude from (\ref{kb19}) that for any fixed $r\in(1,2)$ and $t\in (0,T]$,
\beq\label{kb20}
\|\Delta{\vartheta}\|_{L^r(0,t;L^{\infty})}+\|\Delta{\vartheta}\|_{L^{2r}(0,t;L^2)}\leq C(c_0+|\Delta\vartheta\|_{L^1(0,t;L^{\infty})}+\|\Delta \vartheta\|_{L^2(\Om_t)}).
\eeq
Now the crucial step lies in applying Lemma 3.6 in \cite{ZA1} to (\ref{kb20}) so that we obtain
\beq\label{kb21}
\|\Delta{\vartheta}\|_{L^r(0,t;L^{\infty})}+\|\Delta{\vartheta}\|_{L^{2r}(0,t;L^2)}\leq C c_0.
\eeq
As a consequence, estimate (\ref{f13}) follows readily by combining (\ref{kb19}) and (\ref{kb21}). The proof of Theorem \ref{t2} is thus finished.

\begin{Remark}\label{rem4}
As a matter of fact, like the classical results on one-dimensional compressible Navier-Stokes-Fourier system, the weak solution $(\tau,u,b,\vartheta)$ obtained in Theorem \ref{t1} satisfies
\[
\tau\in C([0,T];L^{\infty}),\,u\in C([0,T];L^2),\,b\in C([0,T];L^{\infty}),\,\vartheta\in C([0,T];L^1),
\]
for any fixed $0<T<\infty$, see \cite{AZ2,ZA1}.
\end{Remark}

\centerline{\bf Acknowledgement}
The research of Yang Li and Yongzhong Sun is supported by NSF of China under Grant No. 11571167 and PAPD of Jiangsu Higher Education Institutions.


\end{document}